\documentclass[12pt]{amsart}
\usepackage{amsmath,amssymb,amsthm,amscd,amsopn,amsfonts,amsxtra} 
\usepackage{latexsym,array}
\usepackage[mathscr]{eucal}
\usepackage{graphicx,psfrag}
\usepackage{epsfig}
\usepackage{color,hyperref}

\newtheorem{theorem}[subsection]{Theorem} 
\newtheorem{lemma}[subsection]{Lemma}

\newtheorem{remark}[subsection]{Remark}

\definecolor{orange}{rgb}{0.995, 0.75, 0.35}
\definecolor{purple}{rgb}{0.7, 0.2, 0.5}
\definecolor{royalblue}{rgb}{0.2, 0.7, 0.8}
\def\al{\alpha}
\def\de{\delta}

\def\lam{\lambda}
\def\om{\omega}

\def\vphi{\varphi}
\def\De{\Delta}
\def\Ga{\Gamma}

\def\Om{\Omega}

\def\th{\tanh}
\def\sech{\mathrm{sech}}

\def\supp{\mathrm{supp}}

\def\inv{^{-1}}
\def\iy{\infty}

\newcommand{\cal}{\mathcal}

\newcommand{\one}{\mathbf{1}}
\newcommand{\la}{\langle}
\newcommand{\ra}{\rangle}
\newcommand{\nd}{\noindent}

\newcommand{\vs}{\vspace}

\newcommand{\hB}{\hfill$\Box$}

\newcommand{\Z}{\mathbb{Z}}
\newcommand{\R}{\mathbb{R}}
\newcommand{\C}{\mathbb{C}}
\newcommand{\N}{\mathbb{N}}

\def\sideremark#1{\ifvmode\leavevmode\fi\vadjust{\vbox to0pt{\vss
 \hbox to 0pt{\hskip\hsize\hskip1em
\vbox{\hsize2cm\tiny\raggedright\pretolerance10000
 \noindent #1\hfill}\hss}\vbox to8pt{\vfil}\vss}}}%

{\end{list}}

\begin{document}
\title[Gradient estimate of heat kernel]
{Note on  gradient estimate of heat kernel for Schr\"odinger operators}

\author{Shijun Zheng}
\address{Department of Mathematical Sciences\\
Georgia Southern University\\
Statesboro, GA 30460-8093}
\email{szheng@georgiasouthern.edu}

\begin{abstract}
Let $H=-\De+V$ be a Schr\"odinger operator on $\R^n$.
We show that  gradient  estimates for the heat kernel of $H$ with upper Gaussian bounds
imply polynomial decay for the kernels of certain smooth dyadic spectral operators. 
The latter decay property has been known to play an important role in the Littlewood-Paley theory for  $L^p$ and Sobolev spaces.  
We are able to  establish the result by modifying Hebisch and the author's recent proofs. 
We give a counterexample in one dimension to show that there exists  $V$ in the Schwartz class such that
 the long time gradient  heat kernel estimate fails.  
\end{abstract}

\maketitle 

\section{Introduction}\label{s:S1}

Consider a Schr\"odinger operator $H=-\De+V$ on $\R^n$, where 
 $V$ is a real-valued potential in $L^1_{loc}(\R^n)$. 
It is noted in \cite{OZ08,Z06a} that for positive $V$,
if $H$ admits the following gradient estimates for its heat kernel $p_t(x,y)=e^{-tH}(x,y)$: 
for all $x,y\in \R^n$ and $t>0$, 
\begin{align}
 &|p_t(x,y)|\le c_n t^{-n/2} e^{-c|x-y|^2/t},  \label{e:pt-gb}\\
& |\nabla_xp_t(x,y)|\le c_n t^{-(n+1)/2}  e^{-c|x-y|^2/t}, \label{e:de-pt-gb}
\end{align}
then the kernel of $\Phi_j(H)$ and its derivatives satisfy a polynomial decay as in (\ref{e:de-phi-j-dec}),
where  $\Phi_j $ is a function in certain Sobolev space with support in $[-2^{j}, 2^{j}]$. 
As is well-known, the decay estimate in (\ref{e:de-phi-j-dec}) implies the Littlewood-Paley  inequality for $L^p(\R^n)$
\cite{E95,D97,OZ06,BZ10}.

For positive $V$, based on heat kernel estimates one can show  
(\ref{e:de-phi-j-dec})
by a scaling argument \cite{Z06a}.  
In this paper we will prove the general case, namely Theorem \ref{t:etH-phi}, by 
modifying the proofs in \cite{He90a, He90b} and \cite{Z06a}.

However, in general 
the gradient estimates (\ref{e:pt-gb}), (\ref{e:de-pt-gb}) do not hold for all $t$.  
This situation may occur when $H$ is a Schr\"odinger operator with negative potential,  
or  the sub-Laplacian on a Lie group of polynomial growth,
cf. \cite{Ou06,FMV06, Gr95,Zh01}. A second part of this paper  is to show such a counterexample,
based on Theorem \ref{t:etH-phi}.

Recall that for a Borel measurable function $\phi: \R\to\C$,   one can define the 
spectral operator $\phi(H)$ by functional calculus
$\phi(H)=\int_{-\infty}^\infty \phi(\lam)dE_\lam   $, where $dE_\lam$ is 
the spectral measure of $H$.  The kernel of $\phi(H)$ is denoted $\phi(H)(x,y)$ in the following sense.
Let $A$ be an operator on a measure space $(M,d\mu)$,  $d\mu$ being a Borel measure on $M$.
 If there exists a locally integrable
function $K_A$: $M\times M\to \C$ such that 
 \begin{align*} \la Af,g\ra=\int_M (Af) gd\mu=\int_{M\times M} K_A (x,y)f(y)g(x)d\mu(x)d\mu(y)
 \end{align*}
for all $f,g$ in $C_0(M)$ with $\supp\, f$ and $\supp\,g$ being disjoint,
where $C_0(M)$ is the set of continuous functions on $M$ with compact supports,
then $A$ is said to have the kernel $A(x,y):=K_A(x,y)$. 
Throughout this paper,  $c$ or $C$ will denote an absolute positive constant. 

The main result is the following theorem for $V\in L^1_{loc}(\R^n)$.
\begin{theorem}\label{t:etH-phi}   
Suppose that the kernel of $e^{-tH}$ satisfies the upper Gaussian bound for $\al=0,1$
\begin{equation}\label{e:der-etH-gb}
|\nabla^\al_x e^{-tH}(x,y)| \le c_n t^{-(n+\al )/2} e^{-c  |x-y|^2/t}\,, \qquad \forall t>0.
\end{equation}
 Let  $\Phi$ be supported in $[-1, 1]$ and belong to
  $H^{\frac{n+1}{2}+N+\de}(\R)$  for some fixed $N\ge 0$ and $\de>0$.
   Then for each $N\ge 0$, there exists a constant $c_{N}$ independent of $\Phi$ such that
for all $j\in\Z$ 
\begin{align}\label{e:de-phi-j-dec}
\vert \nabla_x^\al \Phi_j(H) (x,y) \vert 
\le c_N 2^{j(n+\al )/2} (1+2^{j/2}|x-y|)^{-N}\Vert \Phi\Vert_{H^{\frac{n+1}{2}+N+\de}} \,, 
\end{align}
where  $\Phi_j(x)=\Phi(2^{-j}x)$ and  $H^s:=H^s(\R)$ denotes the usual Sobolev space 
with norm $\Vert f\Vert_{H^s}=\Vert(1-d^2/dx^2)^{s/2}f\Vert_{L^2}$.
\end{theorem} 

When $V$ is positive, a result of the above type was  proved 
and applied to the cases for the Hermite and Laguerre operators \cite{OZ08}. 
The observation was that
 if $V\ge 0$, then the constants corresponding to 
Lemma \ref{l:etH-prop} do not change for $H_\al=-\De+V_\al$
with 
$V_\al(x)= \al^2V(\al x)$ (called {\em scaling-invariance} in what follows),
according to the Feynman-Kac path integral formula
\cite{Si82} \begin{equation*}
e^{-tH}f(x)=E_x\left(f(\om(t))e^{-\int_0^tV(\om(s))ds}\right),
\end{equation*} 
here $E_x $ is the integral over the path space $\Om$ with respect to the Wiener measure $\mu_x$, $x\in\R^n$ 
and $\om(t)$ stands for a brownian motion (generic path).

For general $V$ the technical difficulty is that we do not have such a scaling-invariance. 
We are able to overcome this difficulty by establishing Lemma \ref{l:wei-L1-V},
a scaling version of the weighted $L^1$ inequality
for $\Phi_j(H)(x,y)$ with $\Phi\in H^s$, 
for which we directly use the scaling information 
indicated by the time variable appearing in Lemma \ref{l:etH-prop}. Thus this leads to the proof of 
 the main theorem 
 by combining methods of Hebisch and the author's in \cite{He90a,Z06a}.

In Section 3 we give a counterexample to show
that for $V_\nu(x)=-\nu(\nu+1)(\cosh x)^{-2}$, $\nu\in\N$, the  estimates in (\ref{e:pt-gb}) and (\ref{e:de-pt-gb}) 
fail for  $t\to\iy$.

Note that under the condition in Theorem \ref{t:etH-phi}, (\ref{e:de-phi-j-dec}) is valid for all $\vphi_j\in  C^\infty_0(\R)$, $j\in\Z$ satisfying 
(i) $\supp\; \vphi_j\subset \{ x: |x|\le 2^j\} $ and 
 $\text{(ii)}\; |\vphi_j^{(k)}(x)|\le c_k 2^{-kj}\, ,    \quad \forall j\in \Z$, $k\in \N_0=\{0\}\cup\N$.
A corollary is that 
(\ref{e:der-etH-gb}) implies the Littlewood-Paley inequality
\begin{equation}\label{e:LP-phij}
\Vert f\Vert_{L^p(\R^n)} \approx 
\Vert \big(\sum_{j} |\vphi_j(H)f(\cdot)|^2\big)^{1/2} \Vert_{L^p(\R^n)}\,,\qquad 1<p<\iy
\end{equation}
for both homogeneous and inhomogeneous systems $\{\vphi_j\}$, 
according to \cite[Theorem 1.5]{OZ08}, see also \cite{Z06a,E95}.

\section{heat kernel having upper Gaussian bound implies rapid decay for spectral kernels
}
In this section we  prove  Theorem \ref{t:etH-phi}.
 Following \cite{He90a} we begin with a simple lemma.
\begin{lemma}\label{l:etH-prop} Suppose that (\ref{e:pt-gb}) holds.  
Then \begin{align*} &\int |e^{-tH}(x,y)|^2 dx\le ct^{-n/2}=:\widetilde{C}(t)\\
&\int |e^{-tH}(x,y)| e^{s|x-y|}  dx\le ce^{cs^2t}=: C(s,t).
\end{align*}\end{lemma}

The next lemma  can be easily proved by a duality argument and we omit the details.
\begin{lemma}\label{l:L-2munu}  Let $L$ be a selfadjoint operator on $L^2(\R^n)$
and $\rho,\nu\in L^\iy(\R)$. Then for each $y$,
\begin{align*}\Vert (\rho\nu)(L)(\cdot,y)\Vert_{L^2} \le\Vert \rho(L)\Vert_{2\to 2}\Vert\nu(L)(\cdot,y)\Vert_{L^2}\,.
\end{align*} If in addition $\rho(L)$ is unitary, then the equality holds.
\end{lemma}

Let $w$ be a submultiplicative weight on $\R^n\times\R^n$, i.e., $0\le w(x,y)\le w(x,z)w(z,y)$, $\forall x,y,z\in \R^n$. 
 For simplicity we also assume $w(x,y)=w(y,x)$. Define the norm for  $k\in L^1_{loc}(\R^{2n})$ as follows:
\begin{align*}
\Vert k(x,y)\Vert_w=\sup_{y\in\R^n}\int_{\R^n} |k(x,y)| w(x,y) dx .
\end{align*}
Then given two operators $L_1,L_2$, it holds that
\begin{align}\label{e:L12-w} 
\Vert (L_1L_2)(x,y)\Vert_w\le\Vert L_1(x,y)\Vert_w \Vert L_2(x,y)\Vert_w \,.
\end{align}

The following lemma is a scaling version  of \cite[Theorem 3.1]{He90b} for $V\in L^1_{loc}(\R^n)$. 
\begin{lemma}\label{l:etL-wei-a}  Suppose that (\ref{e:pt-gb}) holds. Let $L_j=e^{-2^{-j}H}$. Then for each $a\ge 0$, 
there exists a constant $c=c(n,a)$ depending on $n,a$ only such that for all $j\in\Z$ and $k\in \R$,
\begin{align}\label{e:eLj-a}
\int \vert (e^{ikL_j}L_j)(x,y)\vert (1+ 2^{j/2}|x-y|)^adx\le c(n,a)(1+|k|)^{n/2+a} .
\end{align}
\end{lemma}
\begin{proof} (a) First we show the case $j=0$. 
For notational convenience write $L=L_0$, 
then  by \mbox{Lemma \ref{l:etH-prop} }
we have, with $t=1$, \begin{align*}
&\sup_y \Vert L(x,y)\Vert_{L^2_x}\le \widetilde{C}(1)=c\\
&\sup_y\int \vert L(x,y) \vert e^{s|x-y|}dx
\le C(s,1)=:C(s)=ce^{cs^2} .
\end{align*}
Let $\phi(x,y)=e^{\beta |x-y|}(1+|x-y|)^a$, $0<\beta<s$.  Then
\begin{align*}
&\Vert L(x,y)\Vert_\phi
= \int \vert L(x,y) \vert e^{s |x-y|} e^{-(s-\beta)|x-y|}(1+|x-y|)^a dx\\
\le& C(s)\sup_{x,y}e^{-(s-\beta)|x-y|}(1+|x-y|)^a=:C(s)c_{s,\beta,a} .
\end{align*}

In view of Lemma \ref{l:L-2munu},  setting  $\ell=\beta\inv |k|\Vert L(x,y)\Vert_\phi\,$, we have 
\begin{align*} 
&\int \vert(e^{ikL}L)(x,y)\vert (1+|x-y|)^a dx=\int_{|x-y|\le \ell}+\int_{|x-y|>\ell}\\
\le&(1+\ell)^{n/2+a}\Vert L(\cdot,y)\Vert_{L^2_x}+\Vert L(x,y)\Vert_\phi e^{-\beta\ell}e^{|k| \Vert L(x,y)\Vert_\phi} \\
\le& \widetilde{C}(1)(1+\beta\inv |k|C(s)c_{s,\beta,a})^{\frac{n}{2}+a}+C(s)c_{s,\beta,a} ,
\end{align*} where we note that by (\ref{e:L12-w})
\begin{align*}
&\Vert e^{ikL}(x,y)\Vert_\phi
\le \sum_{n=0}^\iy \frac{\Vert (ikL)^n(x,y) \Vert_\phi}{n!} \\
\le&\sum_{n=0}^\iy \frac{|k|^n}{n!}\Vert L(x,y)\Vert_\phi^n=e^{|k| \Vert L(x,y)\Vert_\phi} .
\end{align*} 
It is easy to calculate that 
 \begin{align*} c_{s,\beta,a}=\begin{cases} 1& 0\le a\le s-\beta\;\\ 
e^{-(a-s+\beta)}(\frac{a}{s-\beta})^a& a>s-\beta . \quad 
 \end{cases}
 \end{align*}
Hence taking $\beta=s/2$ and fixing $s=s_0>0$ 
give that
\begin{align*} 
&\int \vert (e^{ikL}L)(x,y)\vert (1+ |x-y|)^adx
\le c(s_0,n,a)(1+ |k|)^{\frac{n}{2}+a} .
\end{align*}

(b) Similarly we show the  case for all $j\in\Z$.   
If $L=e^{-2^{-j}H}$, then Lemma \ref{l:etH-prop} tells that with $t=2^{-j}$
\begin{align*} &\int |e^{-2^{-j} H}(x,y)|^2 dx\le c2^{jn/2}\\
&\int |e^{-2^{-j}H}(x,y)| e^{s|x-y|}  dx\le ce^{c 2^{-j}s^2}.
\end{align*}
For $j\in \Z$ let $\phi_j(x,y)=e^{\beta 2^{j/2}|x-y|}(1+2^{j/2}|x-y|)^a$, $0<\beta<s$.  Then similarly to (a)
we obtain \begin{align*}
&\Vert L_j(x,y)\Vert_{\phi_j}=\int |L_j(x,y)| e^{\beta 2^{j/2}|x-y|}(1+2^{j/2}|x-y|)^adx\\
\le& C(s)\sup_{x}e^{-(s-\beta)|x|}(1+|x|)^a
=C(s)c_{s,\beta,a} .
\end{align*}
It follows that, with $\beta=s/2$ and $s=s_0>0$ fixed, 
\begin{align*} 
&\int \vert(e^{ikL_j}L_j)(x,y)\vert (1+2^{j/2}|x-y|)^a dx=\int_{|x-y|\le \ell}+\int_{|x-y|>\ell}\\
\le&c_n\ell^{n/2}(1+2^{j/2}\ell)^{a}\Vert L_j(\cdot,y)\Vert_{L^2_x}+\Vert L_j(x,y)\Vert_{\phi_j} e^{-\beta 2^{j/2}\ell}e^{|k| \Vert L_j(x,y)\Vert_{\phi_j}} \\
\le& c(n,a)(1+|k|)^{n/2+a} ,
\end{align*} where we set $\ell=2^{-j/2}\beta\inv |k|\Vert L_j(x,y)\Vert_{\phi_j}$. 
\end{proof}

We also need a  basic property on the  weighted $\ell^2$ norm of  Fourier coefficients of a compactly supported function in Sobolev space, which can be proved by elementary Fourier expansions.  
\begin{lemma}\label{l:Tg-l2-Hs} Let $s\ge 0$, $T>0$ and $H_0^s([0,T])=\overline{C^\iy_0([0,T])}$ denote the subspace of Sobolev space $H^s(\R)$. Then we have for all $g\in H^s_0([0,T])$, 
 \begin{align}\label{e:T-l2-Hs}
\sqrt{T}\Vert \hat{g}(n)\Vert_{\ell^2_s}\le  c\Vert g\Vert_{H^s_0} \,,
\end{align}
 where $\Vert \{\al_n\} \Vert_{\ell^2_s}=(\sum_{n\in\Z}|\al(n)|^2\la n/T\ra^{2s})^{1/2}$ and 
 $\hat{g}(n)$ are the Fourier coefficents of g over the interval $[0,T]$.  
 \end{lemma}
The inequality in \eqref{e:T-l2-Hs} can be replaced by  equality (however we will not use this improvement), which is a special case of 
the  general norm characterization for periodic  functions in $H^s([0,T])$, see e.g. \cite{SchTr}. 

It follows from Lemma \ref{l:etL-wei-a} and Lemma \ref{l:Tg-l2-Hs} the following weighted $L^1$ estimates for $\Phi_j(H)(x,y)$,  which is an improved version of \cite[Lemma 3.1]{Z06a}, 
where the restriction $V\ge 0$ is removed. 
 \begin{lemma}\label{l:wei-L1-V}  
Suppose $V\in L^1_{loc}(\R^n)$ and the kernel of
 $e^{-tH}$ satisfies  for all $t>0$   
\begin{equation}\label{e:etH-gb}
 |e^{-tH}(x,y)|\le 
 c_n t^{-n/2} e^{- c|x-y|^2/t}\;.
\end{equation}
 If  $s> (n+1)/2+N$,  $N\ge 0$ and 
    $\supp\; \Phi\subset [-10,10]$, then   
 \begin{align*}
\sup_{j\in\Z, \, y\in\R^n}\Vert \Phi(2^{-j}H) (\cdot,y) \la 2^{j/2} (\cdot-y)\ra^N \Vert_{L^1(\R^n)} 
\le c_n \Vert \Phi\Vert_{H^s(\R)} \, ,
\end{align*}
 here $\la x\ra:=1+|x|$. 
\end{lemma} 

\begin{proof} Let $\Phi\subset [-1,1]$. If $\supp\,g\subset I:=[0,2\pi ]$, then $g$ has the Fourier series expansion on $I$
\[ g(x)= \sum_k \hat{g}(k)e^{ikx} ,
\]
where $\hat{g}(k)=\frac{1}{2\pi }\int_0^{2\pi }g(x) e^{-ikx }dx $. 
Let  $\Phi(\lam)=g(e^{-\lam})e^{-\lam}$ and $f_k(\lam)=\lam e^{ik\lam} $. Then $g(y)=\Phi(-\log y)/y$ 
with  $\supp\,g\subset [e\inv,e]$, and so
\begin{align}\label{e:hA-gn-fk-eH} 
\Phi(2^{-j}H)= \sum_k \hat{g}(k)e^{ike^{-2^{-j}H}}e^{-2^{-j}H}=\sum_k \hat{g}(k) {f}_k(e^{-2^{-j}H}) .
\end{align}
It follows from Lemma \ref{l:etL-wei-a}, (\ref{e:hA-gn-fk-eH}) and Lemma \ref{l:Tg-l2-Hs} that for each $y$
\begin{align*} 
&\int |\Phi_j(H)(x,y)|\la 2^{j/2}(x-y)\ra^N dx\le c\sum_k |\hat{g}(k)| (1+|k |)^{n/2+N}\\
=&c\sum_k |\hat{g}(k)| (1+|k |)^{n/2+N+(1+\de)/2} (1+|k |)^{-(\de+1)/2}\\
\le&c\left(\sum_k |\hat{g}(k)|^2 (1+|k |)^{n+2N+1+\de} \right)^{1/2}(\sum_k(1+|k |)^{-\de-1})^{1/2}\\
\le& c \Vert g\Vert_{H_0^{n/2+N+(1+\de)/2}([0,2\pi])}  \de^{-1/2}  \\
\le& c\de^{-1/2}\Vert \Phi\Vert_{H_0^{s}([-1,1])} ,
\end{align*}
where   $\de=s-N-(n+1)/2$ and the last inequality follows from a change of variable and interpolation.
\end{proof} 
 
\begin{remark} Let  $V=V_+-V_-$, $V_\pm\ge 0$ on $\R^n$, $n\ge 3$. 
Then the heat kernel estimate in \eqref{e:etH-gb} holds
if $V_+$ is in Kato class and  $\Vert V_-\Vert_K$,  the global Kato norm of $V_-$, is 
less than $\kappa_n:= \pi^{n/2}/\Ga(\frac{n}{2}-1)$, 
see \cite{DP05}.  Also (\ref{e:etH-gb}) holds whenever $V\ge 0$ is 
locally integrable on $\R^n$,  $n\ge 1$. 
\end{remark}

\subsection{Proof of Theorem \ref{t:etH-phi} }  With (\ref{e:der-etH-gb}) and
Lemma \ref{l:wei-L1-V} we are in a position to prove (\ref{e:de-phi-j-dec}). The proof 
 is similar to that  of Proposition 3.3 in \cite{Z06a} in the case of positive $V$. 
   For completeness, we present the details here.
\begin{align*}
&\nabla^\al_x \Phi_j(H)(x,y)=\int_z \nabla_x^\al e^{-tH}(x,z) (e^{tH}\Phi_j(H) )(z,y) dz.
\end{align*} 
By (\ref{e:der-etH-gb}) we have 
\begin{align*}
&|\nabla^\al_x \Phi_j(H)(x,y)|\\
\le c_n&  t^{-(n+\al )/2} \int e^{-c|x-z|^2/t} \la(x-z)/\sqrt{t}\ra^N 
\la(x-z)/\sqrt{t}\ra^{-N} \la(z-y)/\sqrt{t}\ra^{-N} \\
&\qquad \qquad \quad \cdot \la(z-y)/\sqrt{t}\ra^{N} |(e^{tH}\Phi_j(H) )(z,y) | dz\\ 
\le c_n& t^{-(n+\al )/2} \la(x-y)/\sqrt{t}\ra^{-N} \int \la(z-y)/\sqrt{t}\ra^{N} | (e^{tH}\Phi(2^{-j}H) )(z,y)| dz. 
\end{align*}
Applying Lemma \ref{l:wei-L1-V} with $t=2^{-j}$, 
we obtain 
\begin{align*}
&|\nabla^\al_x \Phi_j(H)(x,y)|\\
\le& c_n t^{-(n+\al )/2} \la(x-y)/\sqrt{t}\ra^{-N}\Vert e^{\lam}\Phi(\lam)\Vert_{H^{\frac{n+1}{2}+N+\de}}\\
\le&c_nt^{-(n+\al )/2} \la(x-y)/\sqrt{t}\ra^{-N}\Vert \Phi\Vert_{H^{\frac{n+1}{2}+N+\de}} \,, \qquad \de>0.
\end{align*}  \hB

\begin{remark}   
In the following section we will show that there exists  $V\in \cal{S}$, the Schwartz class, such that 
(\ref{e:de-phi-j-dec}) does not hold for $j\to \pm\iy$.
 By Theorem \ref{t:etH-phi}, this means that for such $V$ the gradient upper Gaussian bound (\ref{e:der-etH-gb}) 
 does not hold for all $t$. 
\end{remark}

\section{A counterexample to the gradient heat kernel estimate}
Consider the solvable model $H_\nu=-d^2/dx^2 +V_\nu$, $\nu\in\N$, where
\begin{equation*}
V_\nu(x)= -\nu(\nu+1) \sech^2 x.
\end{equation*}
We know from  \cite{OZ06} that solving the Helmholtz equation for $k\in \R\setminus \{0\}$ 
\begin{equation*}
H_\nu e(x,k)=k^2 e(x,k) ,
\end{equation*}   yields the following formula for the continuum eigenfunctions: 
\begin{equation*}
e(x,k)=(\mathrm{sign}(k))^{\nu}\left(\prod_{j=1}^\nu\frac{1}{j+i | k|}
\right)\,P_\nu(x,k)e^{ikx},\end{equation*}
where $P_\nu(x,k)=p_\nu(\th x,ik)$ is defined by the recursion formula
$$ p_\nu (\th x, ik)=
\frac{d}{dx}\big(p_{\nu-1}(\th x, ik)\big)+( ik-\nu\th x) p_{\nu-1}(\th x, ik)\,, 
$$
with $p_0\equiv 1$. 
Note that $e(x,-k)=e(-x,k)$ and  
the function
\begin{equation}\label{e:exey-k}
(x,y,k)\mapsto e(x,k)\overline{e(y,k)}
=\left(\prod_{j=1}^\nu\frac{1}{j^2+ k^2}
\right)\, P_\nu(x,k)P_\nu(y,-k) e^{ik(x-y)} \end{equation}
is real analytic on $\R^3$. 
Moreover, $H_\nu$ has only absolutely continuous spectrum $\sigma_{ac}=[0,\iy)$
and point spectrum $\sigma_{pp}=\{-1,-4,\ldots ,-\nu^2\}\, .$
The corresponding eigenfunctions $\{ e_n\}_{n=1}^\nu$ in $L^2$ are
 Schwartz functions that are linear combinations of 
$\sech^m x \th^\ell x$, $m\in \N$, $\ell\in \N_0$.  

Let $H_{ac}=H_\nu E_{ac}$ 
 denote the absolutely continuous part of $H_\nu$ and $E_{ac}=E_{[0,\iy)}$ the corresponding orthogonal projection.
If $\phi\in C_0(\R)$,  then we have for all $f\in L^1\cap L^2$,    
\begin{equation*}
\phi(H_\nu)f(x)=
\int K(x,y) f(y) dy +  \sum_{n=1 }^\nu \phi(-n^2)(f, e_n) e_n\,,
\end{equation*}
  where $(f, e_n)=\int f(x) \bar{e}_n(x)dx$ and
\begin{equation}\label{eq:kernel-repre}
K(x,y)= (2\pi)^{-1} \int \phi(k^2)e(x,k)\bar{e}(y,k) dk
\end{equation}
 is the kernel of $\phi(H_{ac})=\phi(H) E_{ac}$, cf. \cite{Z04}. 
Since $H_\nu$ has eigenfunctions in $\mathcal{S}(\R)$ 
and $\sigma_{pp}$ is finite,  from now on it is essential to check the kernel 
 $\phi(H_{ac})(x,y)=K(x,y)$  instead  of the kernel of $\phi(H_\nu)$. 

\subsection{Decay for the kernel of $\Phi_j(H)E_{ac}$}  
Let  $\{\vphi_j\}_{j=-\iy}^\iy\subset C^\infty_0(\R)$ satisfy
\quad $\text{(i')}\; 
\supp\; \varphi_j
\subset \{ x: 2^{j-2}\le |x|\le 2^j\} $  and  $\text{(ii')}
 \; |\vphi_j^{(k)}(x)|\le c_k 2^{-kj}\, ,    \forall j\in \Z$, $k\in \N_0$.
    Let $\kappa_j(x,y)=\vphi_j(H_{ac})(x,y)$.
In \cite{OZ06} we showed that  for each $N$ 
\begin{align}\label{e:Kj-dec}
|\kappa_j(x,y)| \le c_N 2^{j/2}(1+2^{j/2}|x-y|)^{-N} , \qquad\forall j\in\Z,
\end{align}   but (with $\al=1$) 
\begin{align}\label{e:der-Kj-dec}
|\partial_x^\al \kappa_j(x,y)| \le c_N 2^{j/2(1+|\al|)}(1+2^{j/2}|x-y|)^{-N}
\end{align}
only holds for $j\ge 0$ and {\em does not} hold for 
all $j<0$.  This suggests that (\ref{e:der-etH-gb}) fails for $\al=1$ and $t>1$ (or more precisely $t\to \iy$), 
according to Theorem \ref{t:etH-phi}.

Now consider the system $\{\Phi_j\}_{j\in\Z}$ which satisfy (i), (ii) as in \mbox{Section \ref{s:S1}.} 
We may assume  $\Phi_j(x)=\Phi(2^{-j}x)$ for a fixed $\Phi$ in $C^\iy([-1,1])$
with $\Phi(x)=1$ on $[-\frac12,\frac12]$.  Let $f^\wedge$ and $f^\vee$ be the Fourier transform and its inverse of $f$ on $\R$. 
 The following lemma shows that  \eqref{e:Kj-dec}  does not hold for $\Phi_j(H_{ac})(x,y)$  when $j\to \iy$. 
\begin{lemma}\label{l:ker-Phi-j-wei} Let  $K_j(x,y)$ be the kernel of $\Phi_j(H)E_{ac}$.
a)  For each $N\in\N_0$ there exists $c_N$ such that for all $ j\le 0$,  
 \begin{align}
\vert  K_j(x,y)\vert \le& c_N 2^{j/2}(1+2^{j/2}|x-y|)^{-N} .
\label{e:Phi-low}  
\end{align}
b) For each $N\in \N_0$  there exists $c_N$ such that for all $ j> 0$,  precisely
\begin{align}
\vert  K_j(x,y)\vert \le& c_N 2^{j/2}(1+2^{j/2}|x-y|)^{-1}  .    
\label{e:Phi-hi}  
\end{align}
In particular, the decay in \eqref{e:Phi-low} does not hold for all $j>0$ with $N>1$. \\
c) There exist positive constants $C$ and $c$  such that for all $ j\in\Z$,   
\begin{align}  &|K_j(x,y)|\le C  |{\Psi_j }^\vee (x-y)|\notag
+  C\int_{-\iy}^\iy |{\Psi_j }^\vee(u)  | e^{- c|x-y-u|} du ,
\end{align}
where $\Psi_j(k)=\Phi_j(k^2) $  and 
it is easily to see that  for each $N$, there exists $c_N$ such that for all $j$ 
\begin{equation*}
|{\Psi_j }^\vee(x-y)| \le c_N  2^{j/2}(1+2^{j/2}|x-y|)^{-N} . \end{equation*}
\end{lemma}
\begin{proof}  (a) 
Let $\lam=2^{-j/2}$.  
By (\ref{eq:kernel-repre}), (\ref{e:exey-k}) and integration by parts
\begin{align*}
&2\pi (i(x-y) )^N K_j(x,y)\\
=&  (-1)^N\int
e^{ik(x-y)}\partial_k^N [\Phi_j(k^2)\prod_{j=1}^\nu({j^2+ k^2})\inv
 P_\nu(x,k)P_\nu(y,-k)]  dk,
\end{align*} 
which can be written as a finite sum of
\begin{equation}\label{e:der-phi-prod-q}
(\th x)^\ell(\th y)^m  \big[ (\Phi_j(k^2))^{(i)}
\big(\prod_{\iota=1}^\nu 
(\iota^2+k^2)^{-1}\big)^{(r)}   (q_{2\nu}(k))^{(s)} \big]^\vee(x-y)
\end{equation}
$0\le \ell, m\le \nu$, 
$i+r+s=N$, $q_{2\nu}(k)$ are polynomials of degree $\le 2\nu$.
We obtain for each $N$ and all $j\le 0$
\begin{align*}
\vert   (x-y)^N K_j(x,y)\vert
= O(\lam^{i-1})=O(\lam^{N-1})= O( 2^{-j/2(N-1)} ) ,\qquad 
\end{align*}
using  \[\begin{cases}
(\Phi_j(k^2))^{(i)}= O(\lam^i)\\
 \big(\prod_{\iota=1}^\nu (\iota^2+k^2)^{-1}\big)^{(r)}=O(\la k\ra^{-2\nu-r})\\
 q_{2\nu}^{(s)}=O(\la k\ra^{2\nu-s}) .
\end{cases}\] This proves (\ref{e:Phi-low}) for $j\le 0$. 

In order to show part (c) for $j\in\Z$, 
using partial fractions we 
 write $K_j(x,y)$ as a finite sum of 
\begin{align}\label{e:ker-Kj-Phi_j}
  &(\th x)^\ell(\th y)^m  \big[ \Phi_j(k^2)
\prod_{\iota=1}^\nu 
(\iota^2+k^2)^{-1} q_{2\nu}(k) \big]^\vee(x-y) ,\end{align}
which is bounded by (up to a constant multiple) \begin{align*}
 |\big[ \Psi_j(k)\big]^\vee(x-y)|  
    +\sum_{\iota=1}^\nu |\big[ \Psi_j(k) \frac{a_{\iota}+b_{\iota}k}{\iota^2+k^2}
  \big]^\vee(x-y)| \,,
\end{align*}  
where $a_\iota, b_\iota\in \R$.
The general term in the  sum is estimated by  
\begin{align*} \big\vert\big[  \Psi_j ( k)
 \frac{a_{\iota} +b_{\iota} k}{\iota^2+k^2 }  \big]^\vee(x-y)\big\vert
\le C  \int |\Psi_j^\vee( u) | e^{- c |x-y-u|} du\,,
 \end{align*}
in terms of the identities 
\begin{align}
&(e^{-|x|})^\wedge (k)=\frac{2}{1+k^2}\label{e:e-fourier}\\
 &(\mathrm{sign}(x) e^{-|x|})^\wedge (k)=\frac{-2ik}{1+k^2}\label{e:sgn-e-fourier} \,.
\end{align} 

 (b)  Finally we prove the sharp estimate in (\ref{e:Phi-hi}).
 For $j>0$, (\ref{e:Phi-low}) does not hold for $N\ge 2$, instead we have only, with $N=0,1$,
\begin{align*}
\vert K_j(x,y)\vert \le c
  2^{j/2}(1+2^{j/2}|x-y|)^{-N} ,  
\end{align*}
by using similar argument and noting 
(\ref{e:der-phi-prod-q}), (\ref{e:e-fourier}), (\ref{e:sgn-e-fourier}).
Indeed, let $J>0$, $N\ge 2$. 
It is easy to find $\{\phi_j\}$ satisfying (i') and (ii') such that 
 \[\Phi_J(x)= 1-\sum_{J}^\infty \phi_j(x) .\]
We have by (\ref{e:Kj-dec})
\[
 \sum_{J}^\infty (x-y)^N |\phi_j(H_{ac})(x,y)| \le c_N\sum_{J}^\infty 2^{-j/2(N-1)}\sim 2^{-J/2(N-1)} .
 \]
On the other hand, from \eqref{e:ker-Kj-Phi_j} and the relevant steps in part (c)  we observe that
if $N>1$, 
\begin{align}\label{e:H_ac-expression}
 &(x-y)^N\one_{[0,\iy)}(H_\nu)(x,y) \notag\\
=&(x-y)^N\int_{-\iy}^\iy e_\nu(x,k)\bar{e}_\nu(y,k)dk \notag\\
=&\, \text{finite sum of } (x-y)^N \sum_\iota \left(\al_\iota e^{-c_\iota|x-y|}+
\mathrm{sign}(x-y)\beta_\iota e^{-c_\iota|x-y|}\right),
\end{align}
where $\al_\iota,\beta_\iota\neq 0$ are of the form $c\th^\ell x\th^m y$.
This shows that the term $ (x-y)^N\Phi_J(H_{ac})(x,y) $  
cannot admit a decay of $2^{-J/2(N-1)}$ for all $J>0$, 
otherwise one would have
\begin{align*}
 &|(x-y)^N\one_{[0,\iy)}(H_\nu)(x,y)|  \lesssim\, 2^{-J/2(N-1)} ,
\end{align*}
 which leads to a contradiction that the sum of those
functions in \eqref{e:H_ac-expression} must vanish,  
  by letting $J\to \iy$.
\end{proof}
\begin{remark} The argument in the proof of part (b) can be made rigorous by replacing $\one_{[0,\iy]}(H_\nu)$ with $\Phi_L(H_\nu)$,
and then let $L\to \iy$ to get the same contradiction. 
\end{remark}

\subsection{The derivative of the kernel of $\Phi_j(H)E_{ac}$}\label{ss:der-ker} Similar argument show that 
\begin{align*}  |\partial_x K_j(x,y)|\le c_N 2^{-j/2(N-2)} / |x-y|^N
\end{align*}
 holds for  all $j>0$ 
 but does not hold for all $j<0$. 
Therefore the inequality in \eqref{e:de-phi-j-dec} does not hold for general $V\in L^1_{loc}$. 

\vs{.23in}
\nd{\bf Acknowledgment}. \quad The author gratefully thanks 
the hospitality and support of Department of  Mathematics, University of California, Riverside 
during his visit in November 2008. In particular he would like to thank Professor Qi S. Zhang for 
the kind invitation and discussions.  The author also wishes  to thank the referee for  careful reading of the original manuscript, whose comments have improved the presentation of this article.

\end{document}